\documentclass[reqno, 12pt, draft]{amsart}
\usepackage{amsmath, amsfonts, amsthm, amssymb}

\numberwithin{equation}{section}

\newtheorem{lemma}{Lemma}[section]
\newtheorem{prop}[lemma]{Proposition}
\newtheorem{thm}[lemma]{Theorem}
\newtheorem{cor}[lemma]{Corollary}
\newtheorem{observation}[lemma]{Observation}

\theoremstyle{definition}

\theoremstyle{remark}

\newcommand{\HH}{\mathbb{H}}
\newcommand{\KK}{\mathbb{K}}
\newcommand{\CC}{\mathbb{C}}
\newcommand{\RR}{\mathbb{R}}
\newcommand{\bbE}{{\mathbf{E}}}
\newcommand{\bbF}{{\mathbf{F}}}
\newcommand{\bbG}{{\mathbf{G}}}
\DeclareMathOperator{\Pf}{Pf}
\DeclareMathOperator{\Alt}{Alt}
\DeclareMathOperator{\Sym}{Sym}
\DeclareMathOperator{\Mat}{M}

\DeclareMathOperator{\sgn}{sgn}
\DeclareMathOperator{\gr}{gr}
\DeclareMathOperator{\Ker}{Ker}

\newcommand{\lie}[1]{\mathfrak{#1}}
\newcommand{\liegl}{\lie{gl}}
\newcommand{\liesp}{\lie{sp}}

\newcommand{\lieo}{\lie{o}}
\newcommand{\lies}{\lie{s}}
\newcommand{\lieg}{\lie{g}}
\newcommand{\liek}{\lie{k}}
\newcommand{\liep}{\lie{p}}
\newcommand{\liepp}{\lie{p}^+}
\newcommand{\liepm}{\lie{p}^-}
\newcommand{\liem}{\lie{m}}
\newcommand{\lieh}{\lie{h}}

\newcommand{\Symm}{{\mathfrak{S}}}
\newcommand{\imgunit}{\sqrt{-1}}

\newcommand{\T}{\,^t\!}
\newcommand{\comb}[2]{{\mathcal{I}^{#1}_{#2}}}
\newcommand{\BAR}[1]{\overline{#1}}

\newcommand{\C}{{\mathbb C}}
\newcommand{\R}{{\mathbb R}}
\newcommand{\Z}{{\mathbb Z}}

\newcommand{\pmtnk}{\C[\rM_{2n,k}]}
\newcommand{\rM}{{\mathrm{M}}}

\newcommand{\Sptnn}{ {\mathrm{Sp}}_{2n}  }
\newcommand{\Glnn}{\GL_n}
\newcommand{\Glkk}{\GL_k}
\newcommand{\Gltnn}{\GL_{2n}}
\newcommand{\Otnn}{{\mathrm{O}}_{2n}}

\newcommand{\hgln}{{\mathcal H}(\Glnn)}
\newcommand{\hotn}{{\mathcal H}(\Otnn)}
\newcommand{\hsptn}{{\mathcal H}(\Sptnn)}

\newcommand{\mathcalh}{{\mathcal H}}

\newcommand{\GL}{{ \mathrm{GL}} }
\newcommand{\hy}{{\hat{y}}}
\newcommand{\hY}{{\hat{Y}}}

\newcommand{\alga}{{\mathfrak{A}}}
\newcommand{\pmkk}{{\C[\rM_{k}]}}

\newcommand{\transpose}[1]{\,{}^t{#1}}
\newcommand{\diag}{\mathop\mathrm{diag}\nolimits{}}
\newcommand{\conjugate}[1]{\overline{\rule{0pt}{1.2ex} #1}}

\newcommand{\U}{\mathrm{U}}
\newcommand{\Sp}{\mathrm{Sp}}
\newcommand{\Mp}{\mathrm{Mp}}
\newcommand{\OO}{\mathrm{O}}
\newcommand{\USp}{\mathrm{USp}}
\newcommand{\composit}{\odot}
\newcommand{\partition}{\mathcal{P}}
\newcommand{\length}[1]{\ell(#1)}

\newcommand{\version}{Ver.~0.0}
\newcommand{\setversion}[1]{\renewcommand{\version}{Ver.~{#1}}}
\setversion{0.1 [Thu Sep  1 19:42:45 JST 2005]}
\setversion{0.2 [Wed Sep  7 17:23:37 JST 2005]}
\setversion{0.3 [Mon Sep 12 11:58:03 JST 2005]}
\setversion{0.4 [Tue Sep 13 15:38:14 JST 2005]}
\setversion{1.0 [Thu Sep 15 12:02:21 JST 2005]}
\setversion{1.1 [Thu Sep 15 18:42:55 JST 2005]}
\setversion{1.2 [Thu Sep 15 18:42:55 JST 2005]}

\title[Harmonics and Capelli identities]
{Intersection of harmonics and Capelli identities for symmetric pairs}

\author[Soo Teck Lee]{Soo Teck Lee$^{\ast}$}
\address{Department of Mathematics, National University of Singapore, 2 Science
Drive 2, Singapore 117543, Singapore.}
\email{matleest@nus.edu.sg}
\thanks{{\tiny$\ast)$ Partially supported by NUS grant R-146-000-067-112.}}

\author[Kyo Nishiyama]{Kyo Nishiyama$^{\ast\ast}$}
\address{
Department of Mathematics\\
\ 
Graduate School of Science\\
\ 
Kyoto University\\
Sa\-kyo, Kyoto 606-8502, Japan}
\email{kyo@math.kyoto-u.ac.jp}
\thanks{{\tiny$\ast\ast)$ Partially supported by JSPS Grant-in-Aid for Scientific Research, {\#}17340037.}}

\author[Akihito Wachi]{Akihito Wachi$^{\ast\ast\ast}$}
\address{Division of Comprehensive Education\\
Hokkaido Institute of Technology\\
Maeda, Teine-ku, Sapporo 006-8585, Japan}
\email{wachi@hit.ac.jp}
\thanks{{\tiny$\ast$$\ast$$\ast$) Partially supported by JSPS Grant-in-Aid for Scientific Research, {\#}16740022.}}


\begin{document}

\maketitle

\begin{abstract}
We consider a see-saw pair consisting of a Hermitian symmetric pair
$(G_\R, K_\R)$ and a compact symmetric pair $(M_\R, H_\R)$, 
where $ (G_\R, H_\R) $ and $ (K_\R, M_\R) $ form 
a real reductive dual pair in a large symplectic group. 
In this setting, 
we get Capelli identities which explicitly represent certain $K_\C$-invariant
elements in $U(\lieg_\C)$ in terms of $H_\C$-invariant elements in $U(\liem_\C)$.
The corresponding $H_\C$-invariant elements are called Capelli elements.

We also give a decomposition of the intersection of
$\OO_{2n}$-harmonics and $\Sp_{2n}$-harmonics as a module of
$\GL_n = \OO_{2n} \cap \Sp_{2n}$, and construct a basis for the
$\GL_n$ highest weight vectors.  This intersection is in the
kernel of our Capelli elements.  
\end{abstract}

\section*{Introduction}

Consider a see-saw pair of real reductive Lie groups 
in the real symplectic group $\Sp_{2N}(\RR)$,
\[
  \begin{array}{ccc}
    G_0  &        & M_0  \\
    \cup & \times & \cup \\
    K_0  &        & H_0,
  \end{array}
\]
where both $(G_0, H_0)$ and $(K_0, M_0)$ are dual pairs (cf.~\cite[(5.1)]{MR90k:22016}).  
Recall that 
the pair $(G_0, H_0)$ is a dual pair 
if $G_0$ and $H_0$ are the mutual commutants of each other in $ \Sp_{2N}(\RR)$.
Let $\lieg_0$ be the Lie algebra of $G_0$
and $\lieg$ its complexification.  
We use the same notation for the other Lie groups.  

The symplectic group $ \Sp_{2N}(\RR) $ has the unique non-trivial double cover $ \Mp_{2N}(\RR) $, 
which is called a metaplectic group.  
It has a unique unitary representation $\omega$ called the Weil representation (or the oscillator representation;  
see \cite{MR90k:22016} or \cite{MR0463359}, for example).  
We use the same symbol $ \omega $ to denote the differentiated action on the Harish-Chandra module of $ \omega $.  
Thus $ \omega $ is a representation of the complexified Lie algebra $\liesp_{2N}$ of $ \Sp_{2N}(\RR)$.  
The Harish-Chandra module can be explicitly realized as the the polynomial ring 
on the $ N $-dimensional vector space $V \simeq \CC^N$ 
(see \cite[p.~537]{MR90k:22016}).  
In this realization, elements of $ \liesp_{2N} $ act by polynomial coefficient differential operators on $ V $.
%
According to the result of Howe (\cite{MR986027}), we have the following equation of the algebras of differential operators:
\begin{equation}
  \label{intro:eq:w(U(g)K)=w(U(m)H)}
  \omega( U(\lieg)^K ) = \omega( U(\liem)^H ),
\end{equation}
where $K$ and $H$ denote the complexifications of $K_0$ and $H_0$ respectively,
and $U(\lieg)^K$ denotes the set of $K$-invariants in the universal enveloping algebra $U(\lieg)$.

Let $\lieg = \liek \oplus \liep$ be the complexified Cartan decomposition.
The subalgebra $S(\liep)^K$ of the $K$-invariants in the symmetric algebra $S(\liep)$ 
is isomorphic to a polynomial ring (see \cite{MR0257285}),
and let $X_1, X_2, \ldots, X_r$ be a set of homogeneous generators of $S(\liep)^K$.
If $\iota: S(\liep) \to U(\lieg)$ is a $K$-linear pseudo-symmetrization map 
(for the definition, see \S~\ref{subsection:pseudo-symmetrization.map}),
then the image $\iota(X_d)$ is $K$-invariant and hence $\omega( \iota(X_d) )$ can be expressed in terms of $\omega( U(\liem)^H )$
thanks to \eqref{intro:eq:w(U(g)K)=w(U(m)H)}:
\begin{equation*}
  \omega( \iota(X_d) ) = \omega( C_d)  \qquad  (C_d \in U(\liem)^H).
\end{equation*}
We call this formula a {\em Capelli identity for a symmetric pair}
and $C_d$ the {\em Capelli element} corresponding to $ X_d $.
In fact, the identity gives a relation between differential operators in various kind of determinantal form, 
which is similar to the original Capelli identity (\cite{capelli90:_sur}).  
However, our Capelli elements are {\it not} central in the enveloping algebra, 
and they depend on the choice of the pseudo-symmetrization map $\iota$.
Despite this, the Capelli elements in this paper 
resemble the original Capelli element (\textit{ibid.}) 
or its variants (see 
\cite{MR1116239}, 
\cite{MR2031449}, 
\cite{nazarov00:_capel},
\cite{molev_nazarov99:_capel_lie},
\cite{MR1445347} and 
\cite{wachi}, 
for example).

We give explicit forms of the Capelli elements for the see-saw pair, where $ G_0/ K_0 $ is an irreducible Hermitian symmetric space, 
and $ M_0 $ is compact.  The main results are summarized in Theorems \ref{thm:caseR}, \ref{thm:caseC} and \ref{thm:caseH} respectively.

Our Capelli identities are closely related to the (pluri-)Harmonic polynomials for compact classical Lie groups.  
In fact, if we consider the orthogonal group $ \OO_{2n}(\R) $ and the unitary symplectic group $ \USp_{2n} $ in $ \U_{2n} $, 
it turns out that the harmonic polynomials are annihilated by Capelli elements (\S~\ref{subsec:annihilator.of.Capelli.elements}).  
Thus the joint harmonics for $ \OO_{2n} $ and $ \USp_{2n} $ are in the kernel of different kinds of Capelli elements, 
which are written as an element of the enveloping algebra of $ \U_{2n} $ in determinantal form.  

We are interested in the kernel of the Capelli elements.  
However, we have not found a good characterization via representation theory.  
Instead, we study the intersection of harmonics for $ \OO_{2n} $ and $ \Sp_{2n} $.  
Note that the joint harmonics carries a representation of $ \GL_n = \OO_{2n} \cap \Sp_{2n} $.
In Theorem \ref{thm:decomposition.of.U.cap.V}, 
we give a decomposition of the joint harmonics as a $ \GL_n $-module in terms of the Littlewood-Richardson coefficients.  
We also indicate how a basis for the space of $ \GL_n $ highest weight vectors can be obtained in \S \ref{sec:tensor.product.algebra}. 
This is related to the tensor product algebras constructed in \cite{HTW}
(see also \cite{HL}).

As a result of our study of joint harmonics, we can answer the following problem, 
which is of independent interest.  
Let $ \rho_{2n}^{\lambda} $ be an irreducible finite dimensional representation of $ \GL_{2n} $ with the highest weight 
$ \lambda = ( \lambda_1, \lambda_2, \ldots, \lambda_k , 0, \ldots, 0) \; ( k \leq n /2 ) $.  
Let us consider the subrepresentation of $ \OO_{2n} $ (respectively $ \Sp_{2n} $) 
generated by the highest weight vector of $ \rho_{2n}^{\lambda} $, and 
denote it by $ \sigma_{2n}^{\lambda} $ (respectively by $ \tau_{2n}^{\lambda} $).  
Then we can identify the intersection $ \sigma_{2n}^{\lambda} \cap \tau_{2n}^{\lambda} $ as a subspace of joint harmonics, 
and we get a decomposition as a representation of $ \GL_n = \OO_{2n} \cap \Sp_{2n} $:
\begin{equation*}
\sigma_{2n}^{\lambda} \cap \tau_{2n}^{\lambda} \simeq 
\bigoplus\nolimits_{\mu, \nu \in \partition_k} c_{\mu, \nu}^{\lambda} \, \rho_n^{\mu \composit \nu} , 
\end{equation*}
where $ c_{\mu, \nu}^{\lambda} $ is the Littlewood-Richardson coefficients, 
$ \partition_k $ is the set of partitions with length $ \leq k $,  
and $ \mu \composit \nu $ is given in \eqref{eqn:definition.of.composit} below.
{From} this, we conclude that the intersection is in the kernel of Capelli elements, which is not easy to see a priori.

\subsection*{Acknowledgment}
The first named author expresses his sincere gratitude to
Kyoto University for warm hospitality during his visit in the summer of 2005.

\subsection*{Notation for finite dimensional irreducible representations}

Let us fix the notation for irreducible finite dimensional representations of the complex matrix groups 
$ \GL_n , \OO_{2n} $ and $ \Sp_{2n} $.

Let $ \lambda = ( \lambda_1, \ldots , \lambda_n ) \in \Z^n $ be a dominant integral weight for $ \GL_n $, i.e., 
$ \lambda_1 \geq \lambda_2 \geq \cdots \geq \lambda_n $.  
We denote by $ \rho_n^{\lambda} $ the irreducible finite dimensional representation of $ \GL_n $ with highest weight $ \lambda $.  
Similarly, if $ \lambda_n \geq 0 $, we denote by $ \sigma_{2n}^{\lambda} $ (respectively $ \tau_{2n}^{\lambda} $) 
the irreducible finite dimensional representation of $ \OO_{2n} $ (respectively $ \Sp_{2n} $) with highest weight $ \lambda $.
We define $ \lambda^{\ast} = - (\lambda_n, \lambda_{n - 1} , \ldots, \lambda_1 ) $, which is the highest weight for 
the contragredient representation $ (\rho_n^{\lambda})^{\ast} $.  
Let 
\begin{equation}
\partition_n = \{ \lambda = ( \lambda_1, \ldots, \lambda_n ) \in \Z_{ \geq 0}^n \mid 
                          \lambda_1 \geq \lambda_2 \geq \cdots \geq \lambda_n \geq 0 \} 
\end{equation}
be the set of partitions of length $ \leq n $.  
We denote the length of a partition $ \lambda $ by $ \length{\lambda} $.  
Take two partitions $ \mu \in \partition_p $ and $ \nu \in \partition_q \; ( p + q \leq n ) $.  Then we denote 
\begin{equation}
\label{eqn:definition.of.composit}
\begin{split}
\mu \composit \nu &= ( \mu , 0, \ldots, 0 , \nu^{\ast} ) \\
&= ( \mu_1 , \mu_2 , \ldots, \mu_p , 0, \ldots, 0, - \nu_q, - \nu_{q - 1}, \ldots, - \nu_1 ) \in \Z^n , 
\end{split}
\end{equation}
which is a dominant integral weight for $ \GL_n $.  

\section{Capelli identities for symmetric pairs}
\label{sec:Capelli-id-for-sym-sp}

\subsection{Pseudo-symmetrization map}
\label{subsection:pseudo-symmetrization.map}

Let $ G_0 $ be a real reductive Lie group, and $ K_0 $ its maximal compact subgroup.  
We denote by 
$ \lieg $ and $ \liek $ the complexified Lie algebra of $ G_0 $ and $ K_0 $ respectively.

Let $ U(\lieg) $ be the enveloping algebra of $ \lieg $, and 
denote by $ F_i U(\lieg) $ the standard filtration of $ U(\lieg) $.  
For the symmetric algebra $ S(\lieg) $, we have the direct sum decomposition by the standard grading 
$ S(\lieg) = \bigoplus_{i \geq 0} S_i(\lieg) $.  
Then there is a natural grading map 
$\gr_i: F_i U(\lieg) \to S_i(\lieg)$ with the kernel $ F_{i -1} U(\lieg) $.

Let $\lieg = \liek \oplus \liep$ be the complexified Cartan decomposition.
If a $K$-linear map $\iota: S(\liep) \to U(\lieg)$ satisfies  
$\gr_i( \iota(u) ) = u$ 
for every homogeneous element $u \in S_i(\liep)$,
we call $\iota$ a \emph{pseudo-symmetrization map}.

There are two basic examples of pseudo-symmetrization maps.
First, the full symmetrization map
\begin{equation*}
Y_1 Y_2 \cdots Y_d \mapsto
(1/d!) \sum_{\sigma\in\Symm_d}
Y_{\sigma(1)} Y_{\sigma(2)} \cdots Y_{\sigma(d)} \qquad 
(Y_i \in \liep) .
\end{equation*}  
Second, 
when $(\lieg, \liek)$ is of Hermitian symmetric type,
we have the irreducible decomposition 
$\liep = \liepp \oplus \liepm$ as $K$-modules.
In this case,
\begin{align}
  \label{eq:def-of-iota}
  \begin{array}{rcl}
    \iota: S(\liep) &\to& U(\lieg) \\
    u_1 u_2 &\mapsto& u_1 u_2
    \qquad (u_1\in U(\liepp), u_2\in U(\liepm))
  \end{array}
\end{align}
is a pseudo-symmetrization map,
since both $\liepp$ and $\liepm$ are $ K $-stable abelian Lie algebras.


\subsection{Abstract Capelli identity for see-saw pairs}

As in the introduction, consider a see-saw pair of real reductive Lie groups 
in the real symplectic group $\Sp_{2N}(\RR)$,
\[
  \begin{array}{ccc}
    G_0  &        & M_0  \\
    \cup & \times & \cup \\
    K_0  &        & H_0,
  \end{array}
\]
where both $(G_0, H_0)$ and $(K_0, M_0)$ are dual pairs.    
Let $\omega$ be the Harish-Chandra module of the Weil representation, realized 
as polynomial coefficient differential operators on the $ N $-dimensional vector space $V \simeq \CC^N$.
Due to Howe, we have the following identity of algebras of polynomial differential operators:
\begin{equation}
  \label{eq:w(U(g)K)=w(U(m)H)}
  \omega( U(\lieg)^K ) = \omega( U(\liem)^H ),
\end{equation}
where $K$ and $H$ denote the complexifications of $K_0$ and $H_0$
respectively.  

Let $\lieg = \liek \oplus \liep$ be the complexified Cartan decomposition.
The subalgebra $S(\liep)^K$ of the $K$-invariants in the symmetric algebra $S(\liep)$ is isomorphic to a polynomial ring,
and let $X_1, X_2, \ldots, X_r$ be a set of homogeneous generators of $S(\liep)^K$.
Choose a $ K $-equivariant pseudo-symmetrization map  $\iota: S(\liep) \to U(\lieg)$.  
Then the image $\iota(X_d)$ is $K$-invariant and hence $\omega( \iota(X_d) )$ can be expressed 
in terms of $\omega( U(\liem)^H )$ thanks to Equation \eqref{eq:w(U(g)K)=w(U(m)H)}:
\begin{equation}
  \omega( \iota(X_d) ) = \omega( C_d)  \qquad  (C_d \in U(\liem)^H).
\end{equation}
We call this formula a \emph{Capelli identity for a symmetric pair}
and $C_d$ the \emph{Capelli element} corresponding to $ X_d $.
Note that our Capelli elements are \emph{not} central  
and they depend on the choice of the pseudo-symmetrization map $\iota$.
Even if we fix the map $ \iota $, 
the Capelli element $ C_d $ is \emph{not} uniquely determined in general 
because the Weil representation $ \omega $ has a non-trivial kernel.  
However, if $ M $ is relatively small by comparison with $ K $, it is uniquely determined.

\section{Capelli identities for symmetric pairs of Hermitian type}
\label{sec:Capelli.identities.for.Hermitian.pairs}

Let us assume that $(\lieg, \liek)$ is of Hermitian symmetric type, 
and take a pseudo-symmetrization map $ \iota $ as in \eqref{eq:def-of-iota}.
There are three of such (irreducible) see-saw pairs given in Table \ref{table:3cases} (see \cite{MR90k:22016}).  
Note that $(M_0, H_0)$ is also a symmetric pair in all the three cases.
\begin{table}[ht]
\label{table:3cases}
\caption{see-saw pairs with $G_0$ Hermitian type, $M_0$ compact}
$
  \begin{array}{cccccc}
    & \Sp_{2N}(\RR) & G_0 & K_0 & M_0 & H_0 \\[3pt]
    \hline
    \\[-5pt]
    \text{Case $\RR$} 
    & \Sp_{2k(p+q)}(\RR) & \Sp_{2k}(\RR) & \U_k & \U_n & \OO_n \\[3pt]
    \text{Case $\CC$} 
    & \Sp_{2(p+q)(r+s)}(\RR) & \U_{p, q} & \U_p \times \U_q & \U_n \times \U_n & \U_n \\[3pt]
    \text{Case $\HH$} 
    & \Sp_{4k(p+q)}(\RR) & \OO^*_{2k} & \U_k & \U_{2n} & \USp_n
  \end{array}
$
\end{table}

In the subsequent subsections, we will give explicit form of the Capelli elements for these three cases.  
Thus we have identities of differential operators in the expression of various minor determinants, which is of independent interest.  

We note that we can give the Capelli identities 
also for the cases where $M_0$ is non-compact
by using the Fourier transform, although we do not describe them in this note.  
We thank Hiroyuki Ochiai and Jiro Sekiguchi for pointing out it to us.

\subsection{Case $\RR$}

In this subsection,
we give the Capelli identity for the symmetric pair of Case $\RR$.
We first fix the notation and describe the generators of $S(\liep)^K$ to state the main theorem.  
%
%
For Case $ \RR $, 
a complex Lie algebra $\lieg$ and 
its subalgebras $\liek$ and $\liep^\pm$ are explicitly given as follows.
\begin{align*}
  \lieg &= \liesp_{2k}
  = \left\{
    \begin{pmatrix} H & G \\ F & -\T H \end{pmatrix}
    \;\Big|\;
    \begin{array}{l}
      H \in \liegl_k, \\
      G, F \in \Sym_k
    \end{array}
  \right\},
  &
  \liepp &= \left\{
    \begin{pmatrix} 0&G \\ 0&0 \end{pmatrix}
    \in \lieg 
  \right\},
  \\
  \liek &= \left\{ 
    \begin{pmatrix} H&0 \\ 0&-\T H \end{pmatrix}
    \in \lieg 
  \right\}
  \simeq \liegl_k,
  &
  \liepm &= \left\{
    \begin{pmatrix} 0&0 \\ F&0 \end{pmatrix} 
    \in \lieg 
  \right\},
\end{align*}
\begin{align*}
  H_{ij} &= E_{ij} - E_{k+j, k+i} \in \liek,
  &
  G_{ij} &= E_{i,k+j} + E_{j, k+i} \in \liepp,
  &
  F_{ij} &= E_{k+i,j} + E_{k+j, i} \in \liepm,
\end{align*}
where $ E_{ij} $ denotes the matrix unit with $ 1 $ at its $ (i,j) $-th position and $ 0 $ elsewhere, 
and $ \Sym_k = \Sym_k(\CC)$ the set of the complex symmetric $k \times k$ matrices.
The complex Lie algebra $\liem$ and its subalgebra $\lieh$ are given by
\begin{align*}
  \liem &= \liegl_n, &
  \lieh &= \lieo_n = \{ X \in \liegl_n \;;\; \T X + X = 0_n \}.
\end{align*}
Let 
$V=\Mat_{n,k}(\CC)$ 
be the space of $ n \times k $ matrices 
and denote the linear coordinate functions on $V$ 
and the corresponding differential operators by
\begin{equation*}
  x_{si}, \quad \partial_{si} \qquad (1\le s \le n, 1\le i\le k) .
\end{equation*}

Let $\lies = \liesp_{2kn}$ be the complex symplectic Lie algebra,
in which both $(\lieg, \lieh)$ and $(\liek, \liem)$ form dual pairs.
We have (the Harish-Chandra module of) the Weil representation $\omega$ of $\lies$
on the space $\CC[V]$ of polynomial functions on $V$.  
The explicit representation operators of $\lieg$ and $\liem$ are given as follows:
\begin{align}
  \label{eq:caseR-def-of-omega}
  \begin{array}{lp{2em}l}
    \displaystyle
    \omega(G_{ij}) = \imgunit \sum_{s=1}^n x_{si} x_{sj},
    &&
    \displaystyle
    \omega(F_{ij}) = \imgunit \sum_{s=1}^n \partial_{si}\partial_{sj},
    \\
    \displaystyle
    \omega(H_{ij}) = \sum_{s=1}^n x_{si} \partial_{sj} 
    + \dfrac{n}{2}\delta_{ij},
    &&
    \displaystyle
    \omega(E_{st}) = \sum_{i=1}^k x_{si} \partial_{ti}
    + \dfrac{k}{2}\delta_{st}. 
  \end{array}
\end{align}
\indent
We now recall the structure of $S(\liep)^K$.
Since $\lieg_0$ is of Hermitian type,
$K \simeq \GL_k(\CC)$ acts multiplicity-freely 
both on the symmetric algebra $S(\liepp)$ and on $S(\liepm)$.  
\begin{equation*}
  S(\liepp) = \bigoplus\nolimits_{\mu \in \partition_k} W_\mu, \quad
  S(\liepm) = \bigoplus\nolimits_{\mu \in \partition_k} W^\ast_\mu \; ; \qquad
  W_{\mu} \simeq \rho_k^{\mu} .
\end{equation*}
%
%
Thus we have the expression of $S(\liep)^K$,
\begin{equation*}
  S(\liep)^K = (S(\liepp) \otimes_{\CC} S(\liepm))^K
  = \bigoplus\nolimits_{\mu \in \partition_k} (W_\mu \otimes_{\CC} W^\ast_\mu)^K , 
\end{equation*}
which is isomorphic to a polynomial ring
with $k$ algebraically independent generators.
For $d = 1, 2, \ldots k$,
the $d$-th generator is the basis vector 
of the one-dimensional vector space
$(W_\mu \otimes_{\CC} W^\ast_\mu)^K$
for $\mu = (2,\ldots,2, 0,\ldots,0) = ( 2^d , 0^{k - d}) $.
%
The explicit form of the generators are
\begin{align}
  \label{eq:caseR-generators}
  X^{\RR}_d &=
  \sum_{I,J \in \comb{k}{d}}
  \det \bbG_{IJ} \cdot \det \bbF_{JI} \in S(\liep)^K
  \qquad
  (d = 1, 2, \ldots, k ; \; k = r),
\end{align}
where $\comb{k}{d}$ is the index set defined by
\begin{equation}
\comb{k}{d} = \{ I \subset \{1,2,\ldots,k\} \;|\; \#I = d \} ,
\end{equation}
and $\bbG_{IJ}$ denotes the  $d\times d$ submatrix of
the $k\times k$ matrix $(G_{ij})$
with the rows and the columns chosen from $I$ and $J$ respectively.  
Note that the generators above belong to the symmetric algebra $S(\liep)$,
and that $G_{ij}$ and $F_{i'j'}$ appearing in the generators
commute with each other in this context.

We use $\iota$ defined by (\ref{eq:def-of-iota})
for the pseudo-symmetrization map.
The images $\iota(X^{\RR}_d)$ 
of the generators (\ref{eq:caseR-generators}) of $S(\liep)^K$
look the same as $X^{\RR}_d$ themselves,
except that the images are in $U(\lieg)$.
In the following theorem,
we use the column-determinant for the determinant of a matrix with
non-commutative entries, defined by
\[
  \det(Z_{ij}) =
  \sum_{\sigma\in\Symm_d} \sgn(\sigma) \,
  Z_{\sigma(1)1} Z_{\sigma(2)2} \cdots Z_{\sigma(d)d}.
\]

\begin{thm}
\label{thm:caseR}
For $1 \le d \le \min(k,n)$,
we have the Capelli identities for the symmetric pair of Case $\RR$:
\begin{equation}
\label{eqn:Capelli.identity.of.case.R}
\omega( \iota( X^{\RR}_d )) = \omega( C^{\RR}_d ) , 
\end{equation}
where $C^{\RR}_d \in U(\liem)^H \; (d=1,2,\ldots,n) $
are the Capelli elements defined by sums of products 
of two $d\times d$ minors with entries in $U(\liem)$:
\begin{align}
\label{eq:caseR-def-of-Cd}
\begin{split}
C^{\RR}_d &=
(-1)^d \sum_{S,T \in \comb{n}{d}}
\det( E_{s_a,t_b} + (d-b-1-k/2) \delta_{s_a,t_b} )_{a,b}
\\
& \hspace{9em} \times
\det( E_{s_a,t_b} + (d-b-k/2) \delta_{s_a,t_b} )_{a,b},
\end{split}
\end{align}
where $s_a$ denotes an element of the index set $S$ with
$s_1 < s_2 < \cdots < s_d$, and the similar rule applies to $ t_b $.
%
%
Note that the identity is trivial 
when $n < d \le k$ since the right hand side becomes an empty sum.
\qed
\end{thm}

Note that \eqref{eqn:Capelli.identity.of.case.R} is in fact an identity of differential operators expressed explicitly by the formula \eqref{eq:caseR-def-of-omega}.

To prove Theorem \ref{thm:caseR},
we demonstrate the computation for the commutative principal symbols.  
We first recall a basic lemma.

\begin{lemma}[Cauchy-Binet]
  Let $R$ be a commutative ring and $d \le N$.
  For a $ d \times N $ matrix $A \in \Mat_{d,N}(R)$ with coefficients in $ R $ and $ N \times d $ matrix $B \in \Mat_{N,d}(R)$, we have
  \begin{align*}
    \det AB =
    \sum_{S\in \comb{N}{d}} \det A_{\bullet,S} \det B_{S,\bullet},
  \end{align*}
  where $A_{\bullet, S}$ is the $d \times d$ submatrix of $A$
  in which all the rows are chosen and the columns are chosen by $S$.
  \qed
\end{lemma}

Define $n \times k$ matrices $X$ and $\partial$ by
\begin{align*}
  X        &= (x_{si})_{1\le s \le n, \; 1\le i \le k}, &
  \partial &= (\partial_{si})_{1\le s \le n, \; 1\le i \le k}.
\end{align*}
In the following computation 
\emph{we take the principal symbols},
and we write the principal symbol of $\partial_{si}$
by the same letter $\partial_{si}$.
For $I, J \in \comb{k}{d}$, the lemma above yields
\begin{align*}
  \omega( \det( \bbG_{IJ} ) ) 
  &= 
  \det \left( \textstyle \imgunit \sum\nolimits_{s=1}^n x_{s,i_a} x_{s,j_b}
  \right)_{1 \le a,b \le d}
  \\ &=
  (\imgunit)^d \det ( \T(X_{\bullet,I}) X_{\bullet,J} )
=
  (\imgunit)^d \sum_{S \in \comb{n}{d}}
  \det \T(X_{SI}) \det X_{SJ}.
\end{align*}
Similarly we have 
\begin{align*}
  \omega( \det( \bbF_{JI} ) ) 
  &=
  (\imgunit)^d \sum_{T \in \comb{n}{d}}
  \det \T(\partial_{TJ}) \det \partial_{TI},
\end{align*}
and the equation of matrices
\begin{align*}
\omega( \bbE_{ST}) 
  &=
  (X \T \partial)_{ST}
  \qquad
  (S, T \in \comb{n}{d}),
\end{align*}
where $ \bbE_{ST} = (E_{s_a,t_b})_{1\le a,b \le d} $ is a matrix with coefficients in $\liem$.
Note that the diagonal shift 
appearing in (\ref{eq:caseR-def-of-omega}) vanishes here,
since we are considering only the principal symbols.
Note also that elements in the expressions above 
commute with each other for the same reason, and we have
\begin{align}
  \sum_{I,J \in \comb{k}{d}}
  \omega( \det \bbG_{IJ} \cdot \det \bbF_{JI} )
  &=
  (-1)^d \sum_{I,J} \sum_{S,T \in \comb{n}{d}}
  \det \T(X_{SI}) \det X_{SJ}
  \det \T(\partial_{TJ}) \det \partial_{TI}
  \notag
  \\ &\overset{(\ast)}{=}
  (-1)^d \sum_{I,S,T}
  \det \T(X_{SI}) \det (X \T \partial)_{ST} \det \partial_{TI}
  \notag
  \\ &\overset{(\ast\ast)}{=}
  (-1)^d \sum_{I,S,T}
  \det (X \T \partial)_{ST} \det \T(X_{SI}) \det \partial_{TI}
  \notag
  \\ &\overset{(\ast)}{=}
  (-1)^d \sum_{S,T}
  \det (X \T \partial)_{ST} \det (X \T \partial)_{ST}
  \notag
  \\ &=
  (-1)^d \sum_{S,T}
  \omega( \det \bbE_{ST} \det \bbE_{ST} ).
  \label{eq:caseR-outline-of-proof}
\end{align}
This is nothing but our desired formula in Theorem \ref{thm:caseR}
except that there are no diagonal shifts in the last expression.
The equalities with $(\ast)$ and $(\ast\ast)$ do {\em not} hold 
when we do {\em not} take the principal symbols.  
We prove the non-commutative analogues of these two equalities with diagonal shifts in the following two lemmas.  
So from now on, $ \partial $ denotes the usual differential operator. 

\begin{lemma}
  \label{lemma:caseR-analog-of-*}
  We have the following equation:
  \begin{align*}
    \sum_{J\in\comb{k}{d}} \det X_{SJ} \cdot \det \partial_{TJ}
    &=
    \det \bigl(
      \omega( E_{s_a,t_b} + (d-b-k/2)\delta_{s_a,t_b})
      \bigr)_{1 \le a,b \le d}.
  \end{align*}
\end{lemma}
\begin{proof}
  Since the right-hand side equals
  $ \det (
      \sum_{i=1}^k x_{s_a,i} \partial_{t_b,i}
      + (d-b)\delta_{s_a,t_b} 
  )_{1 \le a,b \le d}$,
  this lemma reduces to the original Capelli identity 
  (\cite{capelli90:_sur}).
\end{proof}

\begin{lemma}
  \label{lemma:caseR-analog-of-**}
  We have the following formula 
  for $u_1, u_2, \ldots, u_d \in \CC$:
\begin{align*}
    &
    \det X_{SI} \cdot
    \det \bigl( 
{\textstyle 
      \sum_{i=1}^k x_{s_a,i} \partial_{t_b,i} + u_b \delta_{s_a,t_b} 
}
      \bigr)_{1 \le a,b \le d}
    \\ &\qquad=
    \det \bigl( 
{\textstyle 
      \sum_{i=1}^k x_{s_a,i} \partial_{t_b,i} +(u_b-1)\delta_{s_a,t_b} 
}
      \bigr)_{1 \le a,b \le d}
    \cdot \det X_{SI}.
\end{align*}
\end{lemma}
\begin{proof}
  For $s \in S$,
  it is easy to see the commutator
  $\left[ \sum_{i=1}^k x_{si} \partial_{ti}, \det X_{SI} \right]$
  is equal to $\delta_{st} \det X_{SI}$.
  Since the determinant is multi-linear (in its column),
  we have the lemma.
\end{proof}

\begin{proof}[Proof of Theorem \ref{thm:caseR}]
{From} the first equality of (\ref{eq:caseR-outline-of-proof}),
we have
\begin{align*}
  \sum_{I,J \in \comb{k}{d}}
  \omega( \det \bbG_{IJ} \cdot \det \bbF_{JI} )
  =
  (-1)^d \sum_{I,J} \sum_{S,T \in \comb{2n}{d}}
  \det \T (X_{SI}) \det X_{SJ} 
  \det \T (\partial_{TJ}) \det \partial_{TI}.
\end{align*}
It follows from Lemma \ref{lemma:caseR-analog-of-*} that
this is equal to
\begin{align*}
  (-1)^d \sum_{I,S,T}
  \det \T (X_{SI})
  \det \left( \omega(E_{s_a,t_b} + (d-b-k/2) \delta_{s_a,t_b})
  \right)_{1\le a,b \le d}
  \det \partial_{TI}.
\end{align*}
By Lemma \ref{lemma:caseR-analog-of-**},
it turns out that the expression above equals
\begin{align*}
  (-1)^d \sum_{I,S,T}
  \det \left( \omega(E_{s_a,t_b} + (d-b-1-k/2) \delta_{s_a,t_b})
  \right)_{1\le a,b \le d}
  \cdot
  \det \T (X_{SI})
  \det \partial_{TI}.
\end{align*}
By using Lemma \ref{lemma:caseR-analog-of-*} again, 
this is equal to
\begin{align*}
  (-1)^d \sum_{S,T}
  \det \left( \omega(E_{s_a,t_b} + (d-b-1-k/2) \delta_{s_a,t_b})
  \right)_{1\le a,b \le d}
  \\ \qquad\qquad
  \times
  \det \left( \omega(E_{s_a,t_b} + (d-b-k/2) \delta_{s_a,t_b})
  \right)_{1\le a,b \le d}.
\end{align*}
We have thus proved Theorem \ref{thm:caseR}.
\end{proof}


We conclude this subsection by proving
the $H$-invariance of the Capelli element $C^{\RR}_d$
defined in (\ref{eq:caseR-def-of-Cd}).

\begin{prop}
  \label{prop:caseR-invariance}
  The Capelli element is $H$-invariant, that is,
  $C^{\RR}_d \in U(\liem)^H$,
  where $H = \OO_n $ is the complex orthogonal group.
\end{prop}
\begin{proof}
Consider an algebra of tensor products  
$W = \bigwedge \CC^n \otimes_{\CC} \bigwedge \CC^n \otimes_{\CC}U(\liem)$, 
which has a natural $ M $-module structure; $ M = \GL_n $ acts on the both $ \CC^n $ 
by the {\em dual} of the natural representation (hence on its exterior product), 
and acts on $ U(\liem) $ by the adjoint action.  
Denote the standard basis of $\CC^n$ 
in the first and the second factor by $ \{ e_t \}_t $ and $ \{ e'_t \}_t $ respectively.
Define $\eta_t(u)$ and $\eta'_t(u)$ in $W$ by
  \begin{align*}
    \eta_t(u) &= {\textstyle \sum_{s=1}^n e_s (E_{st} + u \delta_{st})},
    &
    \eta'_t(u) &= {\textstyle \sum_{s=1}^n e'_s (E_{st} + u \delta_{st})}.
  \end{align*}
Note that we have for $u, v \in \CC$ and $ T \in \comb{n}{d} $,
  \begin{align*}
    \eta_T(u) &= \sum_{S \in \comb{n}{d}} e_S \det \bbE_{ST}(u),
    &
    \eta'_T(v) &= \sum_{S \in \comb{n}{d}} e'_S \det \bbE_{ST}(v) ,
    && 
  \end{align*}
  where 
\begin{equation*}
\eta_T(u) = \eta_{t_1}(u-1) \eta_{t_2}(u-2) \cdots \eta_{t_d}(u-d) , \qquad
e_S = e_{s_1} e_{s_2} \cdots e_{s_d} 
\end{equation*}
(the product is taken in the exterior power) 
  and $\bbE_{ST}(u)$ denotes the $d \times d$ matrix whose $(a,b)$-entry
  is $E_{s_a,t_b} + (u-b) \delta_{s_a,t_b}$.
Let $W_d$ be the $ \GL_n $-submodule of $W$ spanned by 
$ \{ e_T e'_{T'} \mid T, T' \in \comb{n}{d} \} $.
Then it is easy to check that the mapping
\begin{equation*}
      \Delta: W_d  \to  W , \qquad 
      e_T e'_{T'} \mapsto \eta_T(u) \eta'_{T'}(v) \quad 
    (T, T' \in \comb{n}{d})
\end{equation*}
is a $ \GL_n$-homomorphism for $u, v \in \CC$,
hence an $ \OO_n$-homomorphism in particular.
Since the element $\sum_T e_T e'_T \in W_d$ is $ \OO_n$-invariant, 
$ \sum_{T} \eta_T(u) \eta'_{T}(v) $ is also $ \OO_n$-invariant.  

We define the contraction map $\varepsilon$ on $W_d$ by
  $\varepsilon(e_T e'_{T'}) = \delta_{T,T'}$,
  which is naturally extended to $W_d \otimes_{\CC} U(\liem)$.  
Note that $ \varepsilon $ is $ \OO_n$-equivariant.  
Now we conclude that
\begin{align*}
\varepsilon\left( \sum\nolimits_T \eta_T(u) \eta'_T(v) \right) 
&=
    \varepsilon \left(
      \sum\nolimits_{T,S,S'} e_S \det \bbE_{ST}(u) e'_{S'} \det \bbE_{S'T}(v)
    \right)
    \\ &=
    \sum\nolimits_{S,T} \det \bbE_{ST}(u) \det \bbE_{ST}(v)
\end{align*}
is $ \OO_n$-invariant.
\end{proof}

\subsection{Case $\CC$}

In this subsection,
we treat Case $\CC$. 
%
%
We realize Lie algebras $\lieg$, $\liek$ and $\liep^\pm$ as given below, and define 
elements $ H^{(x)}_{ij}, H^{(y)}_{ij}, F_{ij}, G_{ij} $ of these algebras by
\begin{align*}
  \lieg &= \liegl_{p+q}
  = \left\{
    \begin{pmatrix} H^{(x)} & G \\ F & H^{(y)} \end{pmatrix} 
    \;\Big|\;
    \begin{array}{ll}
      H^{(x)} \in \liegl_p, & G \in \Mat_{p,q}(\CC) \\
      H^{(y)} \in \liegl_q, & F \in \Mat_{q,p}(\CC)
    \end{array}
  \right\},
  \\
  \liek &= \left\{ 
    \begin{pmatrix} H^{(x)}&0 \\ 0&H^{(y)} \end{pmatrix}
    \in \lieg 
  \right\}
  \simeq \liegl_p \oplus \liegl_q,
  \\
  \liepp &= \left\{
    \begin{pmatrix} 0&G \\ 0&0 \end{pmatrix}
    \in \lieg 
  \right\},
  \qquad\qquad
  \liepm = \left\{
    \begin{pmatrix} 0&0 \\ F&0 \end{pmatrix} 
    \in \lieg 
  \right\},
\end{align*}
\begin{align*}
  H^{(x)}_{ij} &= E_{ij} \in \liek
  && (1 \le i,j \le p),
  &
  G_{ij} &= E_{i,p+j} \in \liepp 
  & (1 \le i \le p, 1 \le j \le q),
  \\
  H^{(y)}_{ij} &= E_{p+i,p+j} \in \liek
  && (1 \le i,j \le q),
  &
  F_{ij} &= E_{p+i,j} \in \liepm
  & (1 \le i \le q, 1 \le j \le p).
\end{align*}
The Lie algebra $\liem$, its subalgebra $\lieh$,
and elements of $\liem$ are given by
\begin{align*}
  &\liem = \liegl_n \oplus \liegl_n, 
  &
  &\lieh = \{ (X, -\T X) \in \liem \} \simeq \liegl_n, \\
  &E^{(x)}_{st} = (E_{st}, 0) \in \liem,
  &
  &E^{(y)}_{st} = (0, E_{st}) \in \liem
  && (1 \le s,t \le n).
\end{align*}
Set
$V=\Mat_{n,p}(\CC) \oplus \Mat_{n,q}(\CC)$
and denote the linear coordinate functions
on each component of $V$ by
%
%
$x_{si},\; y_{sj}$ \;
$(1\le s \le n, \; 1\le i\le p, \; 1\le j\le q)$,
respectively.
Put $\lies = \liesp_{2(p+q)n}$, 
in which both $(\lieg, \lieh)$ and $(\liek, \liem)$ form dual pairs.
The Weil representation $\omega$ of $\lies$ is realized on $\CC[V]$,
and its explicit forms are given as follows:
\begin{align}
  \label{eq:caseC-def-of-omega}
\begin{array}{llll}
  \omega(H^{(x)}_{ij}) 
  & \displaystyle  \displaystyle = \sum_{s=1}^n x_{si} \frac{\partial}{\partial x_{sj}}
  + \frac{n}{2}\delta_{ij},
  & \displaystyle 
  \omega(H^{(y)}_{ij}) 
  & \displaystyle = -\sum_{s=1}^n y_{sj} \frac{\partial}{\partial y_{si}}
  - \frac{n}{2}\delta_{ij},
  \\
  \omega(G_{ij}) & \displaystyle = \imgunit \sum_{s=1}^n x_{si} y_{sj},
  & \displaystyle 
  \omega(F_{ji}) & \displaystyle = \imgunit \sum_{s=1}^n 
  \frac{\partial}{\partial x_{si}} \frac{\partial}{\partial y_{sj}},
  \\
  \omega(E^{(x)}_{st}) 
  & \displaystyle = \sum_{i=1}^p x_{si} \frac{\partial}{\partial x_{ti}}
  + \frac{p}{2} \delta_{st},
  & \displaystyle 
  \omega(E^{(y)}_{st}) 
  & \displaystyle = \sum_{j=1}^q y_{sj}  \frac{\partial}{\partial y_{tj}}
  + \frac{q}{2} \delta_{st}.
\end{array}
\end{align}
\indent
In the similar notation to Case $\RR$, 
we have the decomposition of $S(\liep)^K$,
\begin{equation*}
  S(\liep)^K =
  \bigoplus\nolimits_{ \mu \in \partition_r } (W_\mu \otimes_{\CC} W^\ast_\mu)^K , \qquad
  W_{\mu} \simeq \rho_p^{\mu} \otimes \rho_q^{\mu}{}^{\ast} , 
\end{equation*}
where 
%
%
$ r = \min(p, q)$.
In fact, $S(\liep)^K$ is isomorphic to a polynomial ring
with $ r $ algebraically independent generators,
and their explicit forms are
\begin{align}
  \label{eq:caseC-generators}
  X^{\CC}_d &=
  \sum_{I \in \comb{p}{d}, \, J \in \comb{q}{d}}
  \det \bbG_{IJ} \cdot \det \bbF_{JI}
  &&
  ( d = 1, 2 \ldots, r; \quad r = \min(p, q) ).
\end{align}
%

\begin{thm}
  \label{thm:caseC}
  For $1 \le d \le \min(p, q, n)$,
  we have the Capelli identities for the symmetric pair of Case $\CC$.
  \begin{align*}
    \omega( \iota( X^{\CC}_d )) = \omega( C^{\CC}_d ),
  \end{align*}
  where $C^{\CC}_d \in U(\liem)^H \; (d=1,2,\ldots,n) $ are
  the Capelli elements defined by sums of products of
  two $d\times d$ minors with entries in $U(\liem)$:
  \begin{align}
    \label{eq:caseC-def-of-Cd}
    \begin{split}
    C^{\CC}_d &=
      (-1)^d \sum_{S,T \in \comb{n}{d}}
      \det( E^{(x)}_{s_a,t_b} + (d-b-p/2) \delta_{s_a,t_b} )_{a,b}
      \\
      & \hspace{8em} \times
      \det( E^{(y)}_{s_a,t_b} + (d-b-q/2) \delta_{s_a,t_b} )_{a,b}.
    \end{split}
  \end{align}
%
%
  Note that the Capelli identity is trivial 
  when $n < d \le \min(p,q)$.  
\end{thm}

\begin{proof}
  As in Case $\RR$, we define the matrices
  \begin{align*}
    X &= (x_{si})_{1 \le s \le n, 1 \le i \le p},
    &
    Y &= (y_{sj})_{1 \le s \le n, 1 \le j \le q},
    \\
    \partial^X &= \left( {\partial}/{\partial x_{si}} \right)
    _{1 \le s \le n, 1 \le i \le p},
    &
    \partial^Y &= \left( {\partial}/{\partial y_{sj}} \right)
    _{1 \le s \le n, 1 \le j \le q},
  \end{align*}
  and we thus have
  \begin{align*}
    \omega(\det( \bbG_{IJ} ))
    &=
    (\imgunit)^d \sum_{S \in \comb{n}{d}} \det \T (X_{SI}) \cdot \det Y_{SJ},
    \\
    \omega(\det( \bbF_{JI} ))
    &=
    (\imgunit)^d \sum_{T \in \comb{n}{d}} 
    \det \T (\partial^Y_{TJ}) \cdot \det \partial^X_{TI},
  \end{align*}
  for $I \in \comb{p}{d}$, $J \in \comb{q}{d}$.
  Using these formulas, we can prove the identity as follows:
  \begin{align*}
    \sum_{I \in \comb{p}{d}, J \in \comb{q}{d}}
    \omega( \det \bbG_{IJ} \cdot & \det \bbF_{JI} ) 
=
    (-1)^d \sum_{I,J} \sum_{S,T \in \comb{n}{d}}
    \det \!{\T(X_{SI})} \, \det Y_{SJ} \,
    \det \!{\T(\partial^Y_{TJ})} \, \det \partial^X_{TI}
    \\ 
&=
    (-1)^d \sum_{I,J,S,T}
    \det X_{SI}  \det \partial^X_{TI} \cdot
    \det Y_{SJ}  \det \partial^Y_{TJ}
    \\ 
&=
    (-1)^d \sum_{S,T \in \comb{n}{d}}
    \omega\bigl(
      \det( E^{(x)}_{s_a,t_b} + (d-b-p/2) \delta_{s_a,t_b} )_{a,b}
      \\[-5pt]
      & \hspace{11em} \times
      \det( E^{(y)}_{s_a,t_b} + (d-b-q/2) \delta_{s_a,t_b} )_{a,b}
      \bigr).
  \end{align*}
The last equality follows from Lemma \ref{lemma:caseR-analog-of-*} by replacing $k$ with $p$ or $q$.
%
The $ H $-invariance of the Capelli element can be proved similarly to 
Proposition \ref{prop:caseR-invariance} of Case $\RR$.
\end{proof}


\subsection{Case $\HH$}

For Case $\HH$, we need Lemma \ref{lemma:ishikawa-wakayama} on the identity of minor Pfaffians 
due to Ishikawa-Wakayama \cite{MR2002m:15011}.
%
%
Define complex Lie algebras $\lieg$, $\liek$ and $\liep^\pm$
and elements $ H_{ij}, F_{ij}, G_{ij} $ of these algebras by
\begin{align*}
  \lieg &= \lieo_{2k}
  = \left\{
    \begin{pmatrix} H & G \\ F & -\T H \end{pmatrix} 
    \;\Big|\;
    H \in \liegl_k, 
    G, F \in \Alt_k
  \right\},
  &
  \liepp &= \left\{
    \begin{pmatrix} 0&G \\ 0&0 \end{pmatrix} 
    \in \lieg 
  \right\},
    \\
  \liek &= \left\{ 
    \begin{pmatrix} H&0 \\ 0&-\T H \end{pmatrix} 
    \in \lieg 
  \right\} 
  \simeq \liegl_k,
&
  \liepm &= \left\{
    \begin{pmatrix} 0&0 \\ F&0 \end{pmatrix} 
    \in \lieg 
  \right\},
\end{align*}
\begin{align*}
  H_{ij} &= E_{ij} - E_{k+j,k+i} \in \liek, 
  \\
  G_{ij} &= E_{i,k+j} - E_{j,k+i} \in \liepp, 
  &
  F_{ij} &= E_{k+i,j} - E_{k+j,i} \in \liepm
  && (1 \le i,j \le k),
\end{align*}
where $ \Alt_k = \Alt_k(\CC)$ denotes 
the set of the complex alternating $k \times k$ matrices.
The complex Lie algebra $\liem$ and its subalgebra $\lieh$ are given by
\begin{align*}
  \liem &= \liegl_{2n},
  &
  \lieh &=
  \left\{
    \begin{pmatrix} H & G \\ F & -\T H \end{pmatrix} 
    \;\Big|\;
    \begin{array}{l}
      H \in \liegl_n, \\
      G, F \in \Sym_n
    \end{array}
  \right\}
  \simeq \liesp_{2n}.
\end{align*}
We put 
$V=\Mat_{2n,k}(\CC)$ 
and denote the linear coordinate functions on $V$
and the corresponding differential operators by
\begin{equation*}
  x_{si}, \quad \partial_{si} \qquad
  (1\le s \le 2n, \; 1\le i\le k),
\end{equation*}
respectively.
Let $\lies = \liesp_{4kn}$
in which both $(\lieg, \lieh)$ and $(\liek, \liem)$ form dual pairs.
The explicit representation operators of the Weil representation $ \omega $ are given as follows.
\begin{align}
  \label{eq:caseH-def-of-omega}
\begin{array}{llll}
  \omega(H_{ij}) 
  & \displaystyle = \sum_{s=1}^{2n} x_{si} \partial_{sj}  + n \delta_{ij},
  & \displaystyle 
  \omega(G_{ij}) & \displaystyle = \imgunit \sum_{s=1}^n 
  (x_{si} x_{\BAR{s}j} - x_{\BAR{s}i} x_{sj}),
  \\
  \omega(F_{ji}) & \displaystyle = \imgunit \sum_{s=1}^n 
  (\partial_{si} \partial_{\BAR{s}j} - 
    \partial_{\BAR{s}i} \partial_{sj}),
  & \displaystyle 
  \omega(E_{st}) 
  & \displaystyle = \sum_{i=1}^k x_{si} \partial_{ti} + \frac{k}{2} \delta_{st},
\end{array}
\end{align}
where $\BAR{s} = s+n$.
The structure of $S(\liep)^K$ is completely similar to Cases $\RR$ and $ \CC $; 
we have the decomposition of $S(\liep)^K$,
\begin{equation*}
  S(\liep)^K =
  \bigoplus\nolimits_{ \mu } (W_\mu \otimes_{\CC} W^\ast_\mu)^K,  \qquad
  W_{\mu} \simeq \rho_k^{\mu} , 
\end{equation*}
where 
%
%
$\mu$ runs 
over the set of all the partitions 
of the form $(\mu_1, \mu_1, \mu_2, \mu_2, \ldots) \in \partition_k $.  
The polynomial ring $S(\liep)^K$ has $ r = \lfloor k/2 \rfloor$ algebraically independent generators
\begin{align}
  \label{eq:caseH-generators}
  X^{\HH}_d &=
  \sum_{I \in \comb{k}{2d}}
  \Pf \bbG_{II} \cdot \Pf \bbF_{II} 
  &&
  ( d = 1, 2, \ldots, r; \quad r =  \lfloor k/2 \rfloor ),
\end{align}
where $\Pf$ denotes the Pfaffian of an alternating matrix.
%

\begin{thm}
  \label{thm:caseH}
  For $1 \le d \le \min(\lfloor k/2 \rfloor, n)$,
  we have the Capelli identities for the symmetric pair of Case $\HH$:
  \begin{align*}
    \omega( \iota( X^{\HH}_d )) = \omega( C^{\HH}_d ),
  \end{align*}
  where $C^{\HH}_d \in U(\liem)^H \; (d=1,2,\ldots,n) $ are
  the Capelli elements defined by sums of 
  $2d \times 2d$ minors with entries in $U(\liem)$:
  \begin{align}
    \label{eq:caseH-def-of-Cd}
    C^{\HH}_d &=
      \sum_{S_0, T_0 \in \comb{n}{d}}
      \det( E_{s_a,t_b} + (2d-b-k/2) \delta_{s_a,t_b} )
      _{1 \le a,b \le 2d},
  \end{align}
  where $S \in \comb{2n}{2d}$ is defined
as $S = \{s_1, s_2, \ldots, s_d; n+s_1, n+s_2, \ldots, n+s_d \}$ in terms of 
$S_0 = \{s_1, s_2, \ldots, s_d \}\in \comb{n}{d}$.
  %
%
  Note that the identity is trivial 
  when $n < d \le \lfloor k/2 \rfloor$.
\end{thm}

To prove the theorem, 
we use the following lemma to compute Pfaffians.
\begin{lemma}[Ishikawa-Wakayama \cite{MR2002m:15011}]
  \label{lemma:ishikawa-wakayama}
  Let $R$ be a commutative ring and $d \le n$.
  For $A,B \in \Mat_{n,2d}(R)$, $X \in \Sym_n(R)$,
  we have
  \begin{align*}
    \Pf(\!\T{A} X B - \!\T{B} X A) = \sum_{S\in\comb{2n}{2d}}
    \Pf {\scriptstyle \begin{pmatrix} 0&X \\ -X&0 \end{pmatrix}}\rule[-1.5ex]{0pt}{4ex}_{SS} 
    \det \text{\small$ \begin{pmatrix} A \\ B \end{pmatrix} $}\rule[-1.5ex]{0pt}{4ex}_{S,\bullet} .
\end{align*}
  In particular, when $X=I_n$ we have
  \begin{align*}
    \Pf( \T AB - \T BA ) = (-1)^{d(d-1)/2} \sum_{S_0\in\comb{n}{d}}
    \det \text{\small$\begin{pmatrix} A \\ B \end{pmatrix}$}\rule[-1.5ex]{0pt}{4ex}_{S,\bullet},
  \end{align*}
  where $S \in \comb{2n}{2d}$ is determined by $ S_0 $ as in Theorem \ref{thm:caseH}.
\end{lemma}

\begin{proof}
  The first formula is due to Ishikawa-Wakayama 
  \cite[Corollary 2.1]{MR2002m:15011}.
  For the second formula we use two facts.
  The condition
  $\Pf \left( \begin{smallmatrix} 0&I_n \\ -I_n&0 \end{smallmatrix}
  \right)_{SS} \ne 0$
  implies $s_{d+i} = n+s_i$ ($1 \le i \le d$),
  and that 
  $\Pf \left( \begin{smallmatrix} 0&I_d \\ -I_d&0 \end{smallmatrix}
  \right) = (-1)^{d(d-1)/2}$.
  These facts give the second formula.
\end{proof}

\begin{proof}[Proof of Theorem \ref{thm:caseH}]
Define matrices by
  \begin{align*}
    X &= (x_{si})_{1 \le s \le n, 1 \le i \le k},
    &
    \BAR{X} &= (x_{\BAR{s}i})_{1 \le s \le n, 1 \le i \le k},
    \\
    \partial &= (\partial_{si})_{1 \le s \le n, 1 \le i \le k},
    &
    \BAR{\partial}
    &= (\partial_{\BAR{s}i})_{1 \le s \le n, 1 \le i \le k}.
  \end{align*}
Then we have the equations of matrices
  \begin{align*}
    \omega(\bbG_{II})
    &=
    \imgunit (\T (X_{\bullet,I}) \BAR{X}_{\bullet,I} -
      \T (\BAR{X}_{\bullet,I}) X_{\bullet,I}),
    \\
    \omega(\bbF_{II})
    &=
    -\imgunit (\T (\partial_{\bullet,I}) \BAR{\partial}_{\bullet,I}
      - \T (\BAR{\partial})_{\bullet,I} \partial_{\bullet,I}),
  \end{align*}
  for $I \in \comb{k}{2d}$.
  Using these equations we prove the identity as follows:
  \begin{align*}
    &
    \sum_{I \in \comb{k}{2d}}
    \omega( \Pf \bbG_{II} \cdot \Pf \bbF_{II} )
    \\ &\quad=
    \sum_{I \in \comb{k}{2d}}
    \Pf (\T (X_{\bullet,I}) \BAR{X}_{\bullet,I} -
      \T (\BAR{X}_{\bullet,I}) X_{\bullet,I})
    \cdot
    \Pf (\T (\partial_{\bullet,I}) \BAR{\partial}_{\bullet,I}
      - \T (\BAR{\partial}_{\bullet,I}) \partial_{\bullet,I})
    \\ &\quad=
    \sum_I \sum_{S_0\in\comb{n}{d}} (-1)^{d(d-1)/2}
    \det \begin{pmatrix} X_{\bullet,I} \\ \BAR{X}_{\bullet,I}
    \end{pmatrix}\rule[-1.5ex]{0pt}{4ex}_{S,\bullet}
    \cdot
    \sum_{T_0 \in \comb{n}{d}} (-1)^{d(d-1)/2}
    \det \begin{pmatrix}
	\partial_{\bullet,I} \\ \BAR{\partial}_{\bullet,I}
    \end{pmatrix}\rule[-1.5ex]{0pt}{4ex}_{T,\bullet}
    \tag{by Lemma~\ref{lemma:ishikawa-wakayama}}
\end{align*}
\begin{align*}
    & \quad=
    \sum_{I,S_0,T_0} 
    \det \text{\small$\begin{pmatrix} X \\ \BAR{X} \end{pmatrix}$}\rule[-1.5ex]{0pt}{4ex}_{SI} \cdot
    \det \text{\small$ \begin{pmatrix}
	\partial \\ \BAR{\partial} 
    \end{pmatrix}$}\rule[-1.5ex]{0pt}{4ex}_{TI}
    \\ &\quad=
    \sum_{S_0, T_0 \in \comb{n}{d}}
    \omega \Bigl( \det(E_{s_a,t_b} + (2d-b-k/2) \delta_{s_a,t_b})
      _{1 \le a,b \le 2d}  \Bigr) .
    \tag{by Lemma \ref{lemma:caseR-analog-of-*}}
  \end{align*}
The $H$-invariance of $C^{\HH}_d$ can be shown similarly as in the proof of Proposition \ref{prop:caseR-invariance}, 
which we omit.
\end{proof}


\section{The intersection of harmonics}
\label{sec:intersection.of.harmonics}

\subsection{Diamond pair}

The three see-saw pairs considered above is a part of the diamond pair listed below.
\begin{equation}
\begin{array}{c@{}c@{}c@{}c@{}c@{}c@{}c@{}c@{}}
 & & \U_{2n} & & \\
 & \diagup & & \diagdown & \\
\USp_{2n} & & & & \OO_{2n}(\R) \\
& \diagdown & & \diagup & \\
 & & \U_{n} & & 
\end{array}
\qquad
\Longleftrightarrow
\qquad
\begin{array}{c@{}c@{}c@{}c@{}c@{}c@{}c@{}c@{}}
 & & \U_{k,k} & & \\
 & \diagup & & \diagdown & \\
\Sp_{2k}(\R) & & & &\OO_{2k}^{\ast} \\
& \diagdown & & \diagup & \\
 & & \U_k & & 
\end{array}
\end{equation}
Note that, in the diagram, each pair which is connected by line is a symmetric pair.   
We give the explicit realization of each group here because the realization itself is important in the sequel.
\begin{align*}
\U_{2n} &= \{ g \in \GL_{2n}(\C) \mid \transpose{\conjugate{g}} = g^{-1} \} \\
\USp_{2n} &= \{ g \in \U_{2n} \mid \transpose{g} J_n g = J_n \} , \quad 
J_n = \text{\footnotesize$\begin{pmatrix} 0 & -1_n \\ 1_n & 0 \end{pmatrix}$} \\
\OO_{2n}(\R) &= \{ g \in \GL_{2n}(\R) \mid \transpose{g} = g^{-1} \} 
          = \U_{2n} \cap \GL_{2n}(\R) \\
\U_n &= \USp_{2n} \cap \OO_{2n}(\R) 
\end{align*}
Note that the unitary group $ \U_n $ above is not realized in the standard form.
\begin{align*}
\U_{k, k} &= \{ g \in \GL_{2k}(\C) \mid \transpose{\conjugate{g}} I_{k, k} g = I_{k, k} \} , \quad
I_{k, k} = \diag ( 1_k , - 1_k ) \\
\Sp_{2k}(\R) &= \{ g \in \U_{k, k} \mid \transpose{g} J_k g = J_k \} \\
\OO_{2k}^{\ast} &= \{ g \in \U_{k, k} \mid \transpose{g} = g^{-1} \} \\
\U_{k} &= \Sp_{2k}(\R) \cap \OO_{2k}^{\ast}
\end{align*}
These are almost standard realizations except $ \U_k $ in the bottom.  
For $ \Sp_{2k}(\R) $, see \cite[\S VII.10, Problem 30]{MR1920389}.  
In the book \cite[\S I.17, Eq.~(1.141)]{MR1920389}, $ \OO_{2k}^{\ast} $ is defined precisely as here.  
Also refer \cite[\S VII.10, Problems 32 \& 35]{MR1920389} for $ \U_k $.

In the following, we complexified the compact subgroups (the left hand side of the diamond pair) and consider the complex groups 
$ \GL_{2n}, \Sp_{2n}, \OO_{2n} $ and $ \GL_n $.

\subsection{A problem}  

Consider an irreducible finite dimensional representation $ \rho_{2n}^{\lambda} $ of $ \GL_{2n} $ 
with the highest weight $ \lambda \in \partition_{2n} $.  
Let $B_{2n} = A_{2n} N_{2n} $ be the 
standard Borel subgroup of upper triangular matrices in $ \GL_{2n} $, where 
$A_{2n}$ is the diagonal torus and 
$N_{2n}$ is the maximal unipotent subgroup consisting of all  the upper triangular matrices with $1$'s on the diagonal.  
We take a nonzero highest weight vector $v_{\lambda}$ in $ \rho_{2n}^{\lambda} $ with respect to $ B_{2n} $.  
Let $ \OO_{2n} $ (respectively $ \Sp_{2n} $) be the orthogonal group (respectively symplectic group) realized as a subgroup of $ \GL_{2n} $ as above.  Let us consider
\begin{equation}
\label{eqn:definition.of.U.V.lambda}
\begin{cases}
U_{\lambda}&=\mbox{irreducible $\Otnn$ submodule in $\rho_{2n}^{\lambda}$ generated by $v_{\lambda}$,}\\
V_{\lambda}&=\mbox{irreducible $\Sptnn$ submodule in $\rho_{2n}^{\lambda}$ generated by $v_{\lambda}$}.
\end{cases}
\end{equation}
Here, if $ \lambda $ is a partition of length $ \leq n $, we have isomorphisms 
\begin{math}
U_{\lambda}\cong\sigma_{2n}^{\lambda} \text{ and } 
V_{\lambda}\cong \tau_{2n}^{\lambda} .
\end{math}

Note that 
\begin{math}
\Otnn\cap \Sptnn \cong \Glnn 
\end{math}
(see Equation \eqref{eqn:intersection.of.O.and.Sp} below).
We shall denote this explicitly realized subgroup of $\Gltnn$ again by $\Glnn$.  
Note that the intersection 
$U_{\lambda}\cap V_{\lambda}$  is a $\Glnn$ module. 
 We shall study the following problems.  
\begin{enumerate}
\item[(i)] What is the $\Glnn$ module structure of $U_{\lambda}\cap V_{\lambda}$?
\item[(ii)] What are the  $\Glnn$ highest weight vectors in $U_{\lambda}\cap V_{\lambda}$?
\item[(iii)] How to characterize the subspace $U_{\lambda}\cap V_{\lambda}$ inside the representation of $ \GL_{2n} $?
\end{enumerate}
The problems (i) and (ii) are closely related to the theory of joint harmonics, and we will give 
a satisfactory answer to these problems.  
On the other hand, the problem (iii) seems to be related to the Capelli elements for symmetric pairs.  
For this, we do not have a complete answer yet, but we will prove that the subspace $ U_{\lambda}\cap V_{\lambda} $ is 
in a joint eigenspace of the Capelli elements.

\subsection{The intersection of harmonics}

{From} now on, we assume $ n \geq k $ throughout the rest of this section.

Choose a joint complete polarization of $ \C^{2n} $ with respect to the standard symmetric and symplectic bilinear forms. 
So we have $ \C^{2n} = L^+ \oplus L^- $, where $ L^{\pm} $ is a maximal totally isotropic space for the both bilinear forms.  
If we choose a basis of $ L^+ $ and then $ L^- $, then they together form a basis for $ \C^{2n} $.  
We can assume that $\Otnn$ is the subgroup of $\Gltnn$ which preserves the symmetric bilinear form 
\begin{equation*}
((x_1,...,x_n,y_1,...,y_n),
(x^\prime_1,...,x^\prime_n,y^\prime_1,...,y^\prime_n))
 =\sum_{j=1}^n(x_jy^\prime_{j}+y_jx^\prime_j)
\end{equation*}
in the coordinate of $\C^{2n}$ with respect to the basis specified above.
Also $\Sptnn$ is the isometry group of the symplectic form $\langle.,.\rangle$ on $ \C^{2n}$  given by
\begin{equation*}
\langle (x_1,...,x_n,y_1,...,y_n),
(x^\prime_1,...,x^\prime_n,y^\prime_1,...,y^\prime_n)\rangle=\sum_{j=1}^n(x_jy^\prime_{j}-y_jx^\prime_j)
\end{equation*}
with respect to the same coordinate of $ \C^{2n} $.  
In the following, we fix this specific basis and use the coordinate expression freely.  
Then, in matrix form, the intersection of $ \OO_{2n} $ and $ \Sp_{2n} $ is given by 
\begin{equation}
\label{eqn:intersection.of.O.and.Sp}
\Otnn\cap \Sptnn=\left\{\left(\begin{array}{cc}
g&0\\
0& \transpose{g}^{-1}
\end{array}\right):\ g\in \Glnn\right\} ,
\end{equation}
which we will denote simply by $ \GL_n $.  
Let $\rM_{2n,k}= \rM_{2n,k}(\C)$ be the space of $2n\times k$ complex matrices, and let 
\begin{equation*}
\begin{pmatrix} X \\ Y \end{pmatrix} \in \rM_{2n,k} ; \quad 
X = ( x_{ij} )_{1 \leq i \leq n, 1 \leq j \leq k} , \quad
Y = ( y_{ij} )_{1 \leq i \leq n, 1 \leq j \leq k} 
\end{equation*}
be the system of standard coordinates on $\rM_{2n,k} $.
For $1\leq i,j\leq n$, put 
\begin{equation}
\label{eq:definition.of.Laplace.operators}
\Delta_{ij} = \sum_{a=1}^n\frac{\partial^2}{\partial x_{ai}\partial y_{aj}} \; , \qquad\qquad 
\begin{cases}
\Delta^O_{ij} = \Delta_{ij}+\Delta_{ji},\\[7pt]
\Delta^{Sp}_{ij} = \Delta_{ij}-\Delta_{ji}.
\end{cases}
\end{equation}
Notice that the Laplace operators above coincide with $\omega(F_{ij})$'s in 
Equations \eqref{eq:caseR-def-of-omega}, \eqref{eq:caseC-def-of-omega} and \eqref{eq:caseH-def-of-omega} up to constant multiple.  
We now define the following spaces of harmonic polynomials:
\begin{eqnarray*}
\hgln&=&\{f\in\pmtnk:\ \Delta_{ij}(f)=0\ \forall 1\leq i,j\leq n\},\\
\hotn&=&\{f\in\pmtnk:\ \Delta^O_{ij}(f)=0\ \forall 1\leq i,j\leq n\},\\
\hsptn&=&\{f\in\pmtnk:\ \Delta^{Sp}_{ij}(f)=0\ \forall 1\leq i,j\leq n\}.
\end{eqnarray*}
It is clear from the definitions that $\hgln=\hotn\cap\hsptn$. 

Let $\Gltnn\times\Glkk$ act on $\rM_{2n,k}$ by
\begin{equation*}
(g,h).T = \transpose{g}^{-1} \cdot T \cdot h^{-1} \; ;
\qquad 
( g, h ) \in \Gltnn \times \Glkk , \; T\in \rM_{2n,k}.
\end{equation*}
This action induces an action of $\Gltnn\times\Glkk$ on 
the algebra $\C[\rM_{2n,k}]$ of polynomial functions on 
$\rM_{2n,k}$.
By the $(\Gltnn,\Glkk)$-duality, we have
\begin{equation*}
\pmtnk \cong \bigoplus\nolimits_{\lambda \in \partition_k }\rho_{2n}^{\lambda} \otimes \rho_k^{\lambda} .
\end{equation*}
Let $B_k = A_k N_k $ be the 
standard Borel subgroup of upper triangular matrices in $ \GL_k $ similarly defined as $ B_{2n} \subset \GL_{2n} $.  
By taking $N_k$-invariants, we obtain
\begin{equation*}
\pmtnk^{N_k} \cong \bigoplus\nolimits_{\lambda \in \partition_k} \rho_{2n}^{\lambda} \otimes \bigl( \rho_k^{\lambda} \bigr)^{N_k}.
\end{equation*}
Note that $\pmtnk^{N_k}$  is  a module for $\Gltnn\times A_k$. 
Let $W_{\lambda}$ denote the the submodule $ \rho_{2n}^{\lambda} \otimes \left( \rho_k^{\lambda} \right)^{N_k} $ 
of $\pmtnk^{N_k}$.  
We denote by $ \psi_k^{\lambda} $ the character of $ A_k $ defined by 
\begin{equation*}
\psi_k^{\lambda}(a) = a^{\lambda} = a_1^{\lambda_1} \cdots a_k^{\lambda_k} \qquad ( a = \diag ( a_1, \ldots, a_k) \in A_k ) .
\end{equation*}
Then $W_{\lambda}$ is the $\psi_k^{\lambda}$-eigenspace for $A_k$ in $\pmtnk$. Since $\dim \bigl( \rho_k^{\lambda} \bigr)^{N_k} = 1$, 
if we ignore the action of $A_k$ and consider only the action of $\Gltnn$ on $W_{\lambda}$, then 
$W_{\lambda}$  is a copy of the representation $\rho_{2n}^{\lambda}$. 

We now  let $\xi_{\lambda}$ be a $\Gltnn\times\Glkk$ joint highest weight vector in
$\pmtnk$ of weight $\psi_{2n}^{\lambda}\times\psi_k^{\lambda}$ (which is unique up to scalar multiples). Then $\xi_{\lambda}$ is a $\Gltnn$ highest weight vector in   $W_{\lambda}$.   
Now $\Otnn$ and $\Sptnn$ act on $\pmtnk$ as subgroups of $\Gltnn$. 
In this realization, we can identify the subspaces $ U_{\lambda} $ and $ V_{\lambda} $ in Equation \eqref{eqn:definition.of.U.V.lambda} with 
\begin{equation*}
\begin{cases}
U_{\lambda}&=\mbox{irreducible $\Otnn$ submodule in $\rho_{2n}^{\lambda}$ generated by $\xi_{\lambda}$},\\
V_{\lambda}&=\mbox{irreducible $\Sptnn$ submodule in $\rho_{2n}^{\lambda}$ generated by $\xi_{\lambda}$}.
\end{cases}
\end{equation*}

The space $\hotn$ of $\Otnn$ harmonic polynomials is a module for $\Otnn\times\Glkk$ with its structure 
given by
\begin{equation*}
\hotn \cong \bigoplus\nolimits_{\mu \in \partition_k}\sigma_{2n}^{\mu}\otimes\rho_k^{\mu}.
\end{equation*}
In addition, 
the joint $\Otnn\times \Glkk$ highest weight vector in 
$\sigma_{2n}^{\mu}\otimes\rho_k^{\mu}$ coincides with the $\Gltnn\times\Glkk$
highest weight vector in $\rho_{2n}^{\mu}\otimes \rho_k^{\mu}$. By taking
$N_k$ invariants, we obtain
\begin{equation*}
\hotn^{N_k} \cong \bigoplus\nolimits_{\mu \in \partition_k} \sigma_{2n}^{\mu}\otimes(\rho_k^{\mu})^{N_k}
\end{equation*}
which is now a module for $\Otnn\times A_k$. From here, we see that the subspace 
$U_{\lambda}$ coincides with the submodule $\sigma_{2n}^{\lambda}\otimes(\rho^F_k)^{N_k}$
of $\hotn$, and it is precisely the   $\psi_k^{\lambda}$-eigenspace 
$\hotn^{N_k}_{\psi_k^{\lambda}}$ of $A_k$ in
$\hotn^{N_k}$.  On the other hand, the space $\hsptn$ is a module for $\Sptnn\times\Glkk$ with its structure 
given by
\begin{equation*}
\hsptn \cong \bigoplus\nolimits_{\mu \in \partition_k} \tau_{2n}^{\mu}\otimes\rho_k^{\mu}.
\end{equation*}
In a similar way as the $\Otnn$ harmonics, the subspace $V_{\lambda}$ is the 
$\psi_k^{\lambda}$-eigenspace $\hsptn^{N_k}_{\psi_k^{\lambda}}$ of $A_k$ in
$\hsptn^{N_k}$.  Consequently, 
\begin{equation}\label{one}
U_{\lambda}\cap V_{\lambda}= \hotn^{N_k}_{\psi_k^{\lambda}}\cap \hsptn^{N_k}_{\psi_k^{\lambda}}.
\end{equation}
Next we consider the space $\hgln$ of $\Glnn$ harmonics. 
It is a $\Glnn\times\Glkk$ module.
Let  $\hgln^{N_k}_{\psi_k^{\lambda}}$ denotes the $\psi_k^{\lambda}$-eigenspace
of $A_k$ in the space $\hgln^{N_k}$ of $N_k$ invariants in $\hgln$. Then
since $\hgln=\hotn\cap\hsptn$,
\begin{equation*}
\hgln^{N_k}_{\psi_k^{\lambda}}=\hotn^{N_k}_{\psi_k^{\lambda}}\cap \hsptn^{N_k}_{\psi_k^{\lambda}}.
\end{equation*}
It follows from this and  Equation (\ref{one}), 
\begin{equation}\label{two}
U_{\lambda}\cap V_{\lambda}=\hgln^{N_k}_{\psi_k^{\lambda}}.
\end{equation}

\subsection{Decomposition of $ U_{\lambda} \cap V_{\lambda} $ as a $ \GL_n $-module}

We now determine the $\Glnn$ module structure of $\hgln^{N_k}_{\psi_k^{\lambda}}$.  
Put $ \partition_{(k, n)} = \{ ( \mu, \nu ) \in \partition_k \times \partition_k \mid \length{\mu} + \length{\nu} \leq n \} $.  
Then, under the action of $\Glnn\times\Glkk$, the $ \GL_n $ harmonics decomposes as follows.
\begin{eqnarray*}
\hgln 
&\cong& \bigoplus\nolimits_{(\mu, \nu) \in \partition_{(k, n)}}  \rho_n^{\mu \composit \nu} \otimes 
        \left(\rho_k^{\mu}\otimes\rho_k^{\nu}\right)\\
&\cong& \bigoplus\nolimits_{(\mu, \nu) \in \partition_{(k, n)}} \rho_n^{\mu \composit \nu} \otimes 
        \left[ \bigoplus\nolimits_{\eta \in \partition_k} c_{\mu, \nu}^{\eta} \, \rho_k^{\eta}\right]\\
&\cong& \bigoplus\nolimits_{\eta \in \partition_k} 
        \left[ \bigoplus\nolimits_{(\mu, \nu) \in \partition_{(k, n)}} 
        c_{\mu, \nu}^{\eta} \, \rho_n^{\mu \composit \nu}\right]\otimes \rho_k^{\eta} .
\end{eqnarray*} 
%
In the above formula, 
we use the Littlewood-Richardson rule and $c_{\mu, \nu}^{\eta}$ denotes the Littlewood-Richardson coefficient 
which counts the number of the Little\-wood-Rich\-ard\-son tableaux 
of shape $\eta/\mu$ and content $\nu$ (see \cite{MR1464693}).  
By taking $N_k$ invariants, we obtain
\begin{equation}\label{hglnnk}
\hgln^{N_k} \cong 
\bigoplus\nolimits_{\eta \in \partition_k} 
        \left[ \bigoplus\nolimits_{(\mu, \nu) \in \partition_{(k, n)}} c_{\mu, \nu}^{\eta} \, \rho_n^{\mu \composit \nu}\right] 
        \otimes (\rho_k^{\eta})^{N_k}.
\end{equation}
It follows that under $\Glnn$, 
\begin{equation*}
 \hgln^{N_k}_{\psi_k^{\lambda}} \cong 
\bigoplus\nolimits_{(\mu, \nu) \in \partition_{(k, n)}} c_{\mu, \nu}^{\lambda} \, \rho_n^{\mu \composit \nu}.
\end{equation*}
Thus we have proved:

\begin{thm}
\label{thm:decomposition.of.U.cap.V}
Assume that $ k \leq n $, and take $ \lambda \in \partition_k $.  
We define the subspaces 
$ U_{\lambda} \simeq \sigma_{2n}^{\lambda} $ and 
$ V_{\lambda} \simeq \tau_{2n}^{\lambda} $ in $ \rho_{2n}^{\lambda} $ 
as in Equation {\upshape\eqref{eqn:definition.of.U.V.lambda}}.  
Then we have a decomposition under the action of $\Glnn$, 
\begin{equation*}
U_{\lambda}\cap V_{\lambda} \cong
\bigoplus\nolimits_{\mu, \nu \in \partition_k} c_{\mu, \nu}^{\lambda} \, \rho_n^{\mu \composit \nu} , 
\end{equation*}
where the summation is taken over the pair $ ( \mu, \nu ) $ such that $ \length{\mu} + \length{\nu} \leq n $.
\end{thm}

\subsection{Annihilator of Capelli elements}
\label{subsec:annihilator.of.Capelli.elements}

In this subsection, we clarify the relation between the intersection of harmonics and Capelli elements 
which we have constructed in \S~\ref{sec:Capelli.identities.for.Hermitian.pairs}.  

Since the harmonics is by definition annihilated by the (generalized) Laplace operators 
\eqref{eq:definition.of.Laplace.operators}, and those Laplace operators are precisely those $ \omega(F_{ij}) $'s, 
we conclude that they are annihilated by the Capelli elements $ \omega(\iota(X_d^{\KK})) = \omega(C_d^{\KK}) $.  
Here, the symbol $ \KK $ stands for $ \R, \C $ or $ \KK $ respectively.  
Note that $ C_d^{\KK} $'s ($ \KK = \R, \C, \HH $) all belong to the enveloping algebra $ U(\liegl_n) $.  
Thus, the intersection of the harmonics $\hgln=\hotn\cap\hsptn$ are annihilated by all the operators 
$ \{ C_d^{\KK} \mid 1 \leq d \leq r_{\KK} ; \; \KK = \RR, \CC, \HH \} $, 
where $ r_{\KK} = k $ for $ \KK = \RR , \CC $ and $ r_{\HH} = \lfloor k/2 \rfloor $.

However, note that $ \omega(C_d^{\KK}) $ does not coincide with the above action of $ \GL_{2n} $ 
on $ \C[ \Mat_{2n, k} ] $ arisen by the left multiplication.
Instead, it differs by the character (of the two fold covering group) of $ \GL_{2n} $ 
which originates from the definition of the Weil representation.  
In the present situation, the character is $ \chi_k = {\det}^{k/2} $ whose highest weight is 
$ d \chi_k = (k/2) ( 1, 1, \ldots, 1) $.  

\begin{prop}
Assume that $ k \leq n $, and 
let $ U_{\lambda} $ and $ V_{\lambda} $ be as in Theorem \ref{thm:decomposition.of.U.cap.V}.
Then the Capelli elements annihilate $ U_{\lambda} \cap V_{\lambda} $ up to the shift of 
the central character $ \chi_k $.  Namely we have 
\begin{equation*}
\widetilde{C}_d^{\KK} (U_{\lambda}\cap V_{\lambda}) = 0 , \quad 
\text{where } \; 
\widetilde{C}_d^{\KK} = d(\rho_{2n}^{\lambda} \otimes \chi_k) (C_d^{\KK}) 
\quad
(1 \leq d \leq r_{\KK} ; \; \KK = \RR, \CC, \HH) .
\end{equation*}
Here, $ d \rho $ stands for the differentiated representation of $ \liegl_{2n} $.
\end{prop}

It is interesting to study the kernel 
$ \bigcap_{d, \KK} \Ker \widetilde{C}_d^{\KK} $ 
in the representation space of $ \rho_{2n}^{\lambda} $.  
However, up to now, no complete characterization is known.

\section{The $\Glkk$ tensor product algebra}
\label{sec:tensor.product.algebra}

In this section, we assume that $ 2 k \leq n $, 
which is stronger than the one in \S~\ref{sec:intersection.of.harmonics}.  
Under the assumption, we shall prove that 
the space $\hgln^{N_n\times N_k}$ of $N_n\times N_k$-invariants in $\hgln$ is an algebra 
isomorphic to the $\Glkk$ tensor product algebra constructed  in \cite{HTW}.  
Using their results, one can obtain  a basis for $\hgln^{N_n\times N_k}$. 

By taking $N_n$-invariants in Equation \eqref{hglnnk}, we obtain  
\begin{equation}\label{hdecom}
\hgln^{N_n\times N_k} \cong 
           \bigoplus\nolimits_{\lambda, \mu, \nu \in \partition_k } c_{\mu, \nu}^{\lambda} \, 
           \left[(\rho_n^{\mu \composit \nu})^{N_n}\otimes(\rho_k^{\lambda})^{N_k}\right].
\end{equation}
Let  $\mathcalh_{(\mu, \nu), \lambda}$ be the 
$\psi_n^{\mu \composit \nu}\times\psi_k^{\lambda}$-eigenspace for $A_n\times A_k$ in $\hgln^{N_n\times N_k}$. 
Then
\begin{equation}\label{hdecom1}
\hgln^{N_n\times N_k} = \bigoplus\nolimits_{\lambda, \mu, \nu \in \partition_k } \mathcalh_{(\mu, \nu), \lambda},
\end{equation}
and by Equation \eqref{hdecom}, $\dim \mathcalh_{(\mu, \nu), \lambda}=c_{\mu, \nu}^{\lambda}$.

For convenience,   we shall replace the coordinates $Y=(y_{ij})$ by $\hY=(\hy_{ij})$, where 
\begin{equation*}
\hy_{ij}=y_{n-i+1,j},\hspace{1in}(1\leq i\leq n,\ 1\leq j\leq k).
\end{equation*}
We now let $\alga_X$, $\alga_\hY$ and $\alga$ be the subalgebras of $\pmtnk$ generated by the following coordinates:
\begin{eqnarray*}
\alga_X&=&\C[x_{ij}:\ 1\leq i,j\leq k] \cong \pmkk\\
\alga_Y&=&\C[\hy_{ij}:\ 1\leq i,j\leq k] \cong \pmkk\\
\alga&=&\C[x_{ij},\ \hy_{ij}:\ 1\leq i,j\leq k] \cong \alga_X\otimes\alga_Y \cong \pmkk\otimes\pmkk , 
\end{eqnarray*}
where $\rM_{k}=\rM_{k}(\C)$ is the space of $k\times k$ complex matrices. By examining the definition of $\hgln$,   we note that 

\begin{observation}
The algebra $\alga$ is contained in the space $\mathcalh(\Glnn)$ of $\Glnn$ harmonics.
\end{observation}

Now two copies of $\Glkk$ are acting on $\rM_{k}$: 
one by left multiplication and the other by right multiplication. Specifically, for $(g,h)\in\Glkk\times\Glkk$ and $T\in\rM_{k}$,
\begin{equation*}
(g,h).T= \transpose{g}^{-1} \cdot T \cdot h^{-1}.
\end{equation*}
 We shall distinguish these two copies of $\Glkk$ by writing them as $L(\Glkk)$ and $R(\Glkk)$ respectively. More generally, if  $H$ is a
subgroup of $\Glkk$, we shall write the corresponding subgroups of $L(\Glkk)$ and $R(\Glkk)$ as $L(H)$ and
$R(H)$ respectively. In particular, $L(N_k)$ is the maximal unipotent subgroup of $L(\Glkk)$. 

\begin{lemma} \label{isom1}
We have the isomorphism
\begin{equation*}
\mathcalh(\Glnn)^{N_n}\cong\alga_X^{L(N_k)}\otimes \alga_Y^{L(N_k)}.
\end{equation*}
In particular, $\mathcalh(\Glnn)^{N_n}$ is a subalgebra of $\pmtnk$.
\end{lemma}

In fact,  
\begin{equation}\label{hndecom}
\mathcalh(\Glnn)^{N_n}
\cong \bigoplus\nolimits_{\mu, \nu} (\rho_n^{\mu \composit \nu})^{N_n} \otimes \rho_k^{\mu} \otimes \rho_k^{\nu}
\end{equation}
as an $A_n\times \Glkk\times\Glkk$ module. 
In particular, the isotypic component 
$(\rho_n^{\mu \composit \nu})^{N_n}\otimes \rho_k^{\mu} \otimes\rho_k^{\nu}$ 
can be characterized as the 
$\psi_n^{\mu \composit \nu}$-eigenspace for $A_n$ in $\mathcalh(\Glnn)^{N_n}$. 
We shall first describe a spanning set for this isotypic component. 

For $1\leq j\leq k$ and $I= \{ i_1, \ldots ,i_j \} \in \comb{n}{j} $ where $ 1 \leq i_1 < \cdots < i_j \leq n $, let 
\begin{align*}
\gamma_j(I) &=
\det X_{J I} ,
&\text{where $ J = \{1,2,\ldots, j \} $ and $ X = ( x_{s, t} )_{1 \leq s \leq n, 1 \leq t \leq k } $} , \\
%
%
\eta_j(I) &= 
\det \widehat{Y}_{J I} ,
&\text{where $ J = \{1,2,\ldots, j \} $ and $ \widehat{Y} = ( \hat{y}_{s, t} )_{1 \leq s \leq n, 1 \leq t \leq k } $} .
\end{align*}
We also write 
\begin{equation*}
\gamma_j=\gamma_j(\{ 1,2,...,j \}) , \qquad \text{and} \qquad 
\eta_j=\eta_j( \{ 1,2,...,j \} ).
\end{equation*}
For $ \mu = (\mu_1,...,\mu_k) \in \partition_k $ and $ \nu =(\nu_1,...,\nu_k) \in \partition_k $, 
define 
\begin{equation}\label{ab}
\begin{array}{lll}
a_j=\mu_{j-1}-\mu_j , \qquad & b_j=\nu_{j-1}-\nu_j \qquad & \text{for $1\leq j\leq k-1$},\\
a_k=\mu_k , \qquad & b_k=\nu_k.
\end{array}
\end{equation}
Then 
$
\xi_{\mu, \nu}=\gamma_1^{a_1}\gamma_2^{a_2}\cdot\cdot\cdot\gamma_k^{a_k}
\eta_1^{b_1}\eta_2^{b_2}\cdot\cdot\cdot
\eta_k^{b_k}
$
is the unique (up to scalar multiple) $\Glnn\times\Glkk\times\Glkk$ highest weight vector in 
$\mathcalh(\Glnn)$ of weight $\psi_n^{\mu \composit \nu}\times\psi_k^{\mu}\times\psi_k^{\nu}$.

\begin{lemma}
\label{spanning}
For $ \mu, \nu \in \partition_k $, the set of polynomials of the form 
\begin{equation}
\label{generator}
\textstyle
\prod\limits_{j=1}^k \left[\gamma_j(I_j)\right]^{a_j} \, \prod\limits_{m=1}^k \left[\eta_m(J_m)\right]^{b_m}
\end{equation}
spans the isotypic component 
$(\rho_n^{\mu \composit \nu})^{N_n}\otimes \rho_k^{\mu} \otimes\rho_k^{\nu}$ in $\mathcalh(\Glnn)$. 
Here $I_j \in \comb{k}{j} $ and $J_m \in \comb{k}{m} $ are subsets of $\{1,2,...,k\}$ with $j$ and 
$m$ elements respectively, and $a_j$ and $b_m$ are given in 
Equation \eqref{ab}
\end{lemma}

\begin{proof}
First we note that all elements of the form \eqref{generator} are $\Glnn$ highest weight vectors of weight $\psi_n^{\mu \composit \nu}$.  
Let $W$ be the subspace of $\mathcalh(\Glnn)$ spanned by them.  
Then $ W $ is contained in the isotypic component $(\rho_n^{\mu \composit \nu})^{N_n}\otimes \rho_k^{\mu}\otimes\rho_k^{\nu}$ 
in $\mathcalh(\Glnn)$.
%
  
On the other hand, we note that the operators which arise from the infinitesimal action of $\Glkk\times\Glkk$ are derivations 
and affect only the second index $j$ of the variables $x_{ij}$ and $y_{ij}$. 
Thus $W$ is invariant under $\Glkk\times\Glkk$. 
Since the space 
$(\rho_n^{\mu \composit \nu})^{N_n}\otimes \rho_k^{\mu} \otimes\rho_k^{\nu}$ 
is an irreducible $\Glkk\times\Glkk$ module, we have
$
W= (\rho_n^{\mu \composit \nu})^{N_n}\otimes \rho_k^{\mu}
\otimes\rho_k^{\nu}.
$
\end{proof}

\begin{cor} 
Assume that $ 2 k \leq n $.  
Then the space $\mathcalh(\Glnn)^{N_n}$ is the subalgebra of $\pmtnk$ generated by all polynomials 
of the form $\gamma_j(I_j)$ and $\eta_m(J_m)$, 
where  
$1\leq j,m\leq k$ and $I_j$ and $J_m$ are subsets of 
$ \{ 1, 2, \ldots, k \} $ with $j$ and $m$ elements respectively.
\end{cor}

\begin{proof}[Proof of Lemma \ref{isom1}.]
By similar reasonings as Lemma \ref{spanning} (see \cite[p.~83]{MR1321638}), 
the algebra $\alga^{L(N_k)}_X$ is generated by the polynomials $\gamma_j(I_j)$ where $1\leq j\leq k$ and 
$I_j\subseteq \{ 1 , \ldots, k \} $ has $j$ elements.  
Similarly $\alga^{L(N_k)}_Y$ is generated by $\eta_m(J_m)$ where $1\leq m\leq k$ and 
$J_m \subseteq \{ 1, \ldots, k \} $ has $m$ elements. 
It follows that the algebra $\alga_X^{L(N_k)}\otimes \alga_Y^{L(N_k)}$ and 
$\mathcalh(\Glnn)^{N_n}$ have the same generators. Hence they are isomorphic. 
\end{proof}

The $\Glkk$ tensor product algebra is defined as 
\begin{equation*}
{\mathrm{TA}}_k := \bigl[ \pmkk^{R(N_k)}\otimes \pmkk^{R(N_k)} \bigr]{\rule{0pt}{2ex}}^{\Delta(L(N_k))} 
\end{equation*}
in \cite{HTW}. A basis for this algebra was also given in the same paper. 

\begin{cor} 
Under the assumption $ 2 k \leq n $, 
we have the isomorphism 
\begin{equation*}
\hgln^{N_n\times N_k} \cong \bigl[ \alga_X^{L(N_k)}\otimes \alga_Y^{L(N_k)} \bigr]{\rule{0pt}{2ex}}^{\Delta{(R(N_k))}}.
\end{equation*}
In particular, $\hgln^{N_n\times N_k} $ is isomorphic to the $\Glkk$ tensor product algebra ${\mathrm{TA}}_k$.
\end{cor}

Note that by switching the roles of $L(\Glkk)$ and $R(\Glkk)$ in ${\mathrm{TA}}_k$, we obtain the algebra in $\hgln^{N_n\times N_k}$. It follows that to construct a basis for $ \hgln^{N_n\times N_k} $, we only need to take the basis constructed in \cite{HTW} and 
replace their coordinates by ours.  We leave the details to the interested readers.



\def\cprime{$'$}
\providecommand{\bysame}{\leavevmode\hbox to3em{\hrulefill}\thinspace}

\end{document}